\documentclass[12pt,draft]{amsart}
\usepackage{amsfonts}
\usepackage{amssymb}

\usepackage[cp1251]{inputenc}
\usepackage[russian,english]{babel}

\newtheorem{teo}{Теорема}[section]
\newtheorem{lem}[teo]{Лемма}

\newtheorem{cor}[teo]{Следствие}
\newtheorem{dfn}[teo]{Определение}
\newtheorem{rk}[teo]{Замечание}

\renewcommand{\le}{\leqslant}
\renewcommand{\ge}{\geqslant}

\def\<{\langle}
\def\>{\rangle}

\def\a{\alpha}
\def\om{\omega}

\def\e{\varepsilon}
\def\f{{\varphi}}
\def\F{{\Phi}}
\def\A{{\mathcal A}}
\def\M{{\mathcal M}}
\def\cN{{\mathcal N}}
 
\begin{document}

\selectlanguage{russian}

\title[Гильбертов $C^*$-модуль с экстремальными свойствами]
{Гильбертов $C^*$-модуль с экстремальными свойствами}

\author{Д.В.Фуфаев}

\thanks{Данная работа была поддержана грантом Фонда развития теоретической физики и
математики «БАЗИС»}

\address{
Московский Центр фундаментальной и прикладной математики, Отделение МГУ,\newline
Механико-математический факультет МГУ имени М.В.Ломоносова}

\email{denis.fufaev@math.msu.ru, fufaevdv@rambler.ru}

\begin{abstract}
Строится пример гильбертова $C^*$-модуля, который показывает, что теорема Троицкого о геометрической сущности $\A$-компактных операторов между гильбертовыми $C^*$-модулями не может быть распространена на случай несчетнопорожденных модулей (даже если рассматривать более сильную равномерную структуру, которая также будет представлена).
Более того, в построенном модуле не существует фреймов.
\end{abstract}

\maketitle

\section*{Введение}

Хорошо известен критерий компактности операторов в гильбертовых пространствах: оператор компактен (т.е. может быть приближен по норме операторами конечного ранга) тогда и только тогда, когда образ единичного шара вполне ограничен. Этот критерий перестает быть верным, если рассматривать гильбертовы $C^*$-модули, т.е. если рассмотреть некоторую $C^*$-алгебру  $\A$ вместо поля комплексных чисел $\mathbb C$ (в этом случае операторы называются $\A$-компактными). Действительно, даже в случае произвольной бесконечномерной унитальной $C^*$-алгебры $\A$, тождественный оператор имеет ранг, равный единице, однако единичный шар не является вполне ограниченным из-за бесконечной размерности. Поэтому встает естественный вопрос: возможно ли описать свойство $\A$-компактности операторов в геометрических терминах?

Первые шаги в этом направлении были сделаны в \cite{KeckicLazovic2018}, \cite{lazovic2018}, но критерий для некоторого широкого класса модулей был получен только в \cite{Troitsky2020JMAA}. 
А именно, было доказано, что если $\cN$ --- счетнопорожденный гильбертов $C^*$-модуль, то допускающий сопряженный оператор $F:\M\to\cN$ является $\A$-компактным тогда и только тогда, когда образ единичного шара $\M$ является вполне ограниченным множеством в $\cN$ относительно определенной равномерной структуры. После чего возник естественный вопрос: может ли этот результат быть обобщен на случай, когда $\cN$ не является счетнопорожденным?

В \cite{TroitFuf2020} было доказано, что $\A$-компактность оператора влечет вполне ограниченность образа единичного шара относительно данной равномерной структуры для произвольного гильбертова $C^*$-модуля $\cN$. Обратное утверждение было доказано для модулей $\cN$ с определенном свойством проективности. Точнее, с тем свойством, что $\cN$ может быть представлен как
ортогональное прямое слагаемое
в
стандартном модуле над унитализацией алгебры, $\dot{\A}$ (которую мы предполагаем равной самой $\A$ в случае, если она унитальна), некоторой мощности --- то есть, в модуле вида $\bigoplus_{\a\in\Lambda}\dot{\A}$. В частности, это верно для модулей, которые могут быть представлены как прямое слагаемое в $\bigoplus_{\a\in\Lambda}{\A}$ в случае, когда $C^*$-алгебра $\A$ счетнопорождена как модуль над собой. Например, когда $\A$ является $\sigma$-унитальной, т.е. имеет счетную аппроксимативную единицу.

Но в случае, когда $\A$ не $\sigma$-унитальна, как будет доказано в настоящей статье, существует пример, для которого данный критерий не выполняется
Точнее, для некоторого класса (даже коммутативных) $C^*$-алгебр $\A$, рассматриваемых как модуль над собой, будет доказано, что тождественный оператор из $\A$ в себя не $\A$-компактный, но единичный шар в $\A$ (и, что то же самое, образ единичного шара под действием тождественного оператора) вполне ограничен относительно рассматриваемой равномерной структуры. Более того, это верно даже для более сильной равномерной структуры, которая обобщает старую естественным образом.

Также, как это было отмечено как в \cite{Troitsky2020JMAA}, так и в \cite{TroitFuf2020}, построение рассматриваемой равномерной структуры тесно связано с теорией фреймов в гильбертовых $C^*$-модулях, которая берет свое начало в работах Франка и Ларсона (\cite{FrankLarson1999}, \cite{FrankLarson2002}), поэтому представляется естественным рассмотреть построенный пример с точки зрения этой теории и показать, что в нем не существует стандартных фреймов. Ранее пример с этим свойством был построен в \cite{HLi2010}.
Конструкция более сильной равномерной структуры, в свою очередь, тесно связана с теорией внешних фреймов, которую
разработали
Арамбашич и Бакич (\cite{ArBak2017}).

В \S \ref{sec:Prelim} мы напоминаем некоторые факты из топологии и теории гильбертовых $C^*$-модулей и устанавливаем некоторые дополнительные свойства. В \S \ref{sec:TotBoundNoAcomp} мы описываем модуль, для которого упомянутый критерий не выполняется --- а именно, вполне ограниченность образа единичного шара не влечет $\A$-компактность. В \S \ref{sec:nofr} доказывается, что данный модуль не имеет фреймов.

Автор выражает благодарность Е.В.Троицкому и А.И.Корчагину за полезные обсуждения.

\section{Предварительные сведения}\label{sec:Prelim}

Напомним некоторые факты о топологических пространствах, которые, скорее всего, хорошо известны специалистам, но
недостаточно хорошо
представлены в литературе

Для локально компактного хаусдорфова топологического пространства $K$ будем обозначать через $\widetilde{K}=K\cup\{t_\infty\}$ его одноточечную компактификацию (см. \cite[29.1]{Munkr}). Мы также часто будем пользоваться тем фактом, что $C^*$-алгебра $C_0(K)$ непрерывных функций, обращающихся в ноль на бесконечности (см. \cite[436I]{Frem4}) изоморфна идеалу в $C^*$-алгебре $C(\widetilde{K})$, состоящему из функций, обращающихся в ноль в точке $t_\infty$.

\begin{dfn}
Топологическое пространство $K$ называется пространством Урысона (или функционально хаусдорфовым), если для любых двух точек $t,y\in K$ найдется такая непрерывная на $K$ функция, что $f(t)=1$ и $f(y)=0$. Каждое пространство Урысона, очевидно, хаусдорфово. Локально компактное хаусдорфово пространство является пространством Урысона, что следует из \cite[3.1.1]{Engel}.
\end{dfn}

Теперь нам нужен следующий результат о продолжениях и сужениях непрерывных функций.

\begin{teo} (Теорема Титце о продолжении) \cite[Теорема 35.1]{Munkr}
Любая непрерывная функция, определенная на замкнутом подмножестве нормального топологического пространства $K$ может быть продолжена до непрерывной функции на $K$. Если функция ограничена, продолженная функция может быть выбрана ограниченной той же постоянной.
\end{teo}

\begin{cor}\label{extrest}
Пусть $K$ --- локально компактное хаусдорфово пространство, $A$ --- замкнутое подмножество $K$ и 
$f\in C_0(A)$,
тогда $f$ может быть продолжена до функции из $C_0(K)$. Кроме того, $f$ ограничена и продолженную функцию можно выбрать ограниченной той же постоянной.
\end{cor}

\begin{proof}
Рассмотрим $\widetilde {K}=K\cup\{t_\infty\}$ --- одноточечную компактификацию $K$. Оно компактно и хаусдорфово, поэтому нормально (\cite[32.3]{Munkr}). Так как $A\cup\{t_\infty\}$ --- замкнутое подмножество $\widetilde {K}$, можно применить теорему Титце и получить непрерывную функцию на $\widetilde {K}$, которая совпадает с $f$ на $A$ и равна нулю на $\{t_\infty\}$, то есть, функцию из $C_0(K)$.
Так как замкнутое подмножество компакта --- компактно, $f$ ограничена.
\end{proof}

\begin{rk}
Если $A$ замкнуто в $\widetilde {K}$, то $A$ и $\{t_\infty\}$ --- непересекающиеся замкнутые множества. Если $A$ не замкнуто в $\widetilde {K}$ и, следовательно, $t_\infty$ --- предельная точка $A$, то полученная функция непрерывна так как $\lim\limits_{A\ni t\to t_\infty}f(t)=0$.
\end{rk}

Начиная с этого момента пусть $K$ --- локально компактное, не $\sigma$-компактное хаусдорфово топологическое пространство, такое что всякое счетное семейство компактных множеств содержится в некотором компактном множестве (другими словами, всякое $\sigma$-компактное подмножество $K$ содержится в компактном множестве в $K$).
В частности, замыкание счетного набора компактных множеств компактно.

\begin{lem}\label{lem:compsup}
Если $f\in C_0(K)$, то найдется компактное множество $K'\subset K$, такое что $f(t)=0$ для всех $t\in K\setminus K'$ (т.е. $f$ имеет компактный носитель).
\end{lem}

\begin{proof}
Для $n\in\mathbb N$ положим $K_n=\{t\in K: |f(t)|\ge\frac{1}{n}\}$.
По определению $C_0(K)$, для всякого $n\in\mathbb N$ найдется компактное множество $Y_n$, такое что $|f(t)|<\frac{1}{n}$ для всех $t\in K\setminus Y_n$. Так как $K_n\subset Y_n$ и $K_n$ замкнуто, $K_n$ тоже компактно. Отсюда следует, что $|f(t)|\ne 0$ только для $t\in\bigcup\limits_{n=1}^\infty K_n$, так что $f(t)=0$ для любого $t\in K\setminus K'$, где $K'$ --- компакт, содержащий все $K_n$.
\end{proof}

Теперь опишем пример топологического пространства с необходимыми нам свойствами.
Напомним некоторые определения, относящиеся к порядковым числам
(см. \cite[главы VI-VII]{KurMos})

\begin{dfn}
1) Порядковый изоморфизм двух упорядоченных множеств $X$, $Y$ --- это биекция $f:X\to Y$, которая сохраняет порядок, т.е. если $x\le y$ в $X$, то $f(x)\le f(y)$ в $Y$.

2) Порядковый тип --- класс эквивалентности упорядоченных множеств по отношению к порядковому изоморфизму.

3) Вполне упорядоченное множество --- это упорядоченное множество, в котором всякое непустое подмножество имеет наименьший элемент.

4) Порядковое число (или ординал) --- это порядковый тип вполне упорядоченного множества.

5) Всякое множество порядковых чисел упорядочено (более того, вполне упорядочено): если $\a$ и $\beta$ --- порядковые типы вполне упорядоченных множеств $A$ и $B$, то $a\le\beta$ тогда и только тогда, когда существует инъекция $f:A\to B$, сохраняющая порядок.

6) Порядковое число $\om_1$ --- это порядковый тип множества всех счетных ординалов (т.е. порядковых типов не более чем счетных множеств). Это множество несчетно. Более того, для любой последовательности счетных ординалов $\{\a_j\}$ найдется счетный ординал $\a$ такой, что $\a_j<\a$.
\end{dfn}

Рассмотрим
$[0,\om_1]$, то есть, множество всех ординалов $\a$ таких, что $\a\le\om_1$. Его можно представить как топологическое пространство с топологией порядка (\cite[3.1.27]{Engel}): база этой топологии состоит из множеств вида $(y,x]=\{z:y<z\le x\}$, где $y<x\le\om_1$, и множества $\{0\}$ 
(где $0$
это порядковый тип пустого множества). Также рассмотрим $[0,\om_1)=[0,\om_1]\setminus\{\om_1\}$.

Пространство $[0,\om_1]$ компактно, в то время как $[0,\om_1)$ локально компактно и не $\sigma$-компактно; они оба хаусдорфовы (\cite[43.4 и 43.10]{Exampl}).

\begin{lem}\label{lem:compcord}
Всякое счетное семейство компактных подмножеств $\{K_n\}$ пространства $[0,\om_1)$ содержится в некотором компакте $K'\subset [0,\om_1)$.
\end{lem}

\begin{proof}
Так как $\om_1$ не является предельной точкой ни для какого $K_n$ и дополнения $G_n=[0,\om_1)\setminus K_n$ открыты, для любого $n\in\mathbb N$ можно найти открытое множество вида $(\a_n,\om_1)\subset G_n$. Точная верхняя грань множества $\{\a_n\}_{n\in\mathbb N}$ --- это порядковое число $\gamma$ меньшее, чем $\om_1$ в силу того, что $[0,\a_n)$ не более чем счетны и счетное объединение не более чем счетных множеств не более чем счетно. Поэтому, $\bigcup\limits_{n=1}^\infty K_n\subset \bigcup\limits_{n=1}^\infty [0,\a_n]\subset [0,\gamma]=K'$ и $(\gamma,\om_1)$ --- непустое открытое множество.
\end{proof}

\begin{rk}
Для $K=[0,\om_1)$ 
лемма \ref{lem:compsup}
была также доказана в \cite[43.12]{Exampl}.
\end{rk}

Напомним кратко основные сведения о гильбертовых $C^*$-модулях и операторах, которые можно найти в 
\cite{Lance},\cite{MTBook},\cite{ManuilovTroit2000JMS}.

\begin{dfn}
\rm
(Правым) предгильбертовым  $C^*$-модулем над $C^*$-алгеброй $\A$
называется $\A$-модуль, снабженный $\A$-\emph{внутренним произведением}
$\<.,.\>:\M\times\M\to \A$, являющимся такой полуторалинейной формой на подлежащем
линейном пространстве, что  для любых $x,y\in\M$, $a\in\A$:
\begin{enumerate}
\item $\<x,x\> \ge 0$;
\item $\<x,x\> = 0$ тогда и только тогда, когда $x=0$;
\item $\<y,x\>=\<x,y\>^*$;
\item $\<x,y\cdot a\>=\<x,y\>a$.
\end{enumerate}

\emph{Гильбертов $C^*$-модуль} --- это такой предгильбертов $C^*$-модуль над $\A$,
который полон по норме $\|x\|=\|\<x,x\>\|^{1/2}$.

Предгильбертов $C^*$-модуль $\M$
называется
 \emph{счетнопорожденным},
если существует такой счетный набор его элементов, что
их $\A$-линейные комбинации плотны  в $\M$.

Мы будем обозначать через $\oplus$ гильбертову сумму гильбертовых
$C^*$-модулей в понятном смысле.
\end{dfn}

\begin{dfn}\label{dfn:operators}
\rm
\emph{Оператор} --- это ограниченный $\A$-гомоморфизм.
Оператор, имеющий сопряженный (в понятном смысле), называется
\emph{допускающим сопряженный} (см. \cite[\S 2.1]{MTBook}).
Обозначим банахово пространство всех операторов
$F: \M\to \cN$ через ${\mathbf{L}}(\M,\cN)$
и банахово пространство операторов, допускающих сопряженный, --- через ${\mathbf{L}}^*(\M,\cN)$.  
\end{dfn}

\begin{dfn}\label{dfn:Acompact}
\rm
Обозначим через
$\theta_{x,y}:\M\to\cN$, где $x\in\cN$ 
и $y\in\M$, \emph{элементарный} $\A$-\emph{компактный} оператор, который
определяется как $\theta_{x,y}(z):=x\<y,z\>$.
Тогда банахово пространство 
$\A$-\emph{компактных} операторов ${\mathbf{K}}(\M,\cN)$ ---
замыкание подпространства, порожденного всеми элементарными $\A$-компактными операторами
в ${\mathbf{L}}(\M,\cN)$.
\end{dfn}

\begin{rk}\label{prop:propAcomp}
\rm
Пространство ${\mathbf{K}}(\M,\cN)$ является инволютивным подпространством 
${\mathbf{L}}^*(\M,\cN)$, $C^*$-алгебра
${\mathbf{K}}(\M,\M)$ является идеалом $C^*$-алгебры ${\mathbf{L}}^*(\M,\M)$,
${\mathbf{K}}(\M,\cN)$ является ${\mathbf{L}}^*(\M,\M)$-${\mathbf{L}}^*(\cN,\cN)$-бимодулем, а также
${\mathbf{K}}(\M,\cN){\mathbf{L}}^*(\cN,\cN)\subseteq {\mathbf{K}}(\M,\cN)$ и
${\mathbf{L}}^*(\M,\M){\mathbf{K}}(\M,\cN)\subseteq {\mathbf{K}}(\M,\cN)$ (см. \cite[\S 2.2]{MTBook}).
\end{rk}

Чтобы определить равномерную стуктуру, а именно, систему полуметрик, которая будет обобщать ранее введенную (\cite[определения 2.4 и 2.9]{Troitsky2020JMAA}), нужно изложить понятия мультипликаторной алгебры и мультипликаторного модуля (подробнее см. \cite{BakGul2003}, \cite{BakGul2004}).

Напомним, что $M(\A)$ --- это $C^*$-алгебра мультипликаторов $C^*$-алгебры $\A$ (см., например, \cite[глава 2]{Lance}), $M(\A)=\A$, если $\A$ унитальна. 
Кроме того, если $K$ --- локально компактное хаусдорфово топологическое пространство и $\A=C_0(K)$,
то $M(\A)=C_b(K)$ --- алгебра всех непрерывных ограниченных функций на $K$.

Для всякого гильбертова $\A$-модуля $\cN$ существует гильбертов $M(\A)$-модуль $M(\cN)$ (который называется мультипликаторным модулем модуля $\cN$), содержащий $\cN$ в качестве идеального подмодуля, ассоциированного с $\A$, т.е. $\cN=M(\cN)\A$. Более того, для любых $x\in\cN$, $y\in M(\cN)$ выполнено $\<x,y\>\in\A$. $M(\cN)=\cN$, если алгебра $\A$ унитальна. 
Также, в силу того, что каждый элемент $x\in\cN$ может быть представлен в виде $y\cdot a$ для некоторых $y\in\cN$, $a\in\A$ (см. \cite[1.3.10]{MTBook}), модуль $\cN$ и всякий его подмодуль могут рассматриваться как $M(\A)$-модули.

Если мы положим $\cN=C_0(K)$ как модуль над собой, то, как и в случае алгебры мультипликаторов, $M(C_0(K))=C_b(K)$.

Представим теперь новую равномерную структуру на подмодулях модуля $\cN$.

\begin{dfn}\label{dfn:admissyst}
\rm
Рассмотрим гильбертов $C^*$-модуль $\cN$ над $\A$. Счетная система его элементов
$X=\{x_{i}\}$ называется (внешней) \emph{допустимой} для
 (возможно, незамкнутого) подмодуля $\cN^0\subseteq  \cN$
(или внешней $\cN^0$-\emph{допустимой}), если
\begin{enumerate}
\item[1)] 
 при каждом $x\in\cN^0$ ряд $\sum_i \<x,x_i\>\<x_i,x\>$ сходится по норме;
\item[2)]
его сумма ограничена $\<x,x\>$, то есть $\sum_i \<x,x_i\>\<x_i,x\> \le \<x,x\>$; 
\item[3)]$\|x_i\|\le 1$ для любого $i$.
\end{enumerate}
\end{dfn}

\begin{rk}
Точно так же, как допустимая система (\cite[определение 2.1]{Troitsky2020JMAA}) представляет собой некоторый аналог ``относительной нормированной последовательности Бесселя'' в духе теории фреймов, внешняя допустимая система может рассматриваться как ``относительная нормированная внешняя последовательность Бесселя'' (см. \cite{ArBak2017}).
\end{rk}

\begin{dfn}
Обозначим через $\F$ счетный набор $\{\f_1,\f_2,\dots\}$
состояний на $\A$. Для каждой пары $(X,\F)$ с внешней $\cN^0$-допустимой $X$
рассмотрим неотрицательную функцию, определенную по формуле
$$
\nu_{X,\F}(x)^2:=\sup_k 
\sum_{i=k}^\infty |\f_k\left(\<x,x_i\>\right)|^2,\quad x\in \cN^0. 
$$
Обозначим множество всех таких функций как
$\mathbb{OSN}(\cN,\cN^0)$. Также будем писать $(X,\F)\in 
\mathbb{OA}(\cN,\cN^0)$
для пар с внешними допустимыми $X$. 
\end{dfn}

\begin{rk}
Возникает естественный вопрос: почему эта функция --- корректно определенная полунорма? Чтобы обосновать это, рассмотрим $\cN^0$ как $M(\A)$-модуль и, применяя \cite[2.3.24]{BratRob}, возьмем $\widehat{\F}=\{\widehat{\f}_1,\widehat{\f}_2,\dots\}$ --- счетный набор состояний на $M(\A)$, продолжающих $\{\f_1,\f_2,\dots\}$.
Ясно, что система $X$ является $\cN^0$-допустимой системой в $M(\A)$-модуле $M(\cN)$. Поэтому функция
$$
\widehat{\nu}_{X,\widehat{\F}}(x)^2:=\sup_k 
\sum_{i=k}^\infty |\widehat{\f}_k\left(\<x,x_i\>\right)|^2,\quad x\in \cN^0. 
$$
является корректно определенной полунормой на $\cN^0$, так как она удовлетворяет условиям \cite[определения 2.1 и 2.4]{Troitsky2020JMAA}, и она обладает всеми соответствующими свойствами. Очевидно, для $x\in\cN^0$ мы имеем $\<x,x_i\>\in\A$ и, следовательно, $\widehat{\f}_k\left(\<x,x_i\>\right)={\f}_k\left(\<x,x_i\>\right)$, откуда $\widehat{\nu}_{X,\widehat{\F}}(x)={\nu}_{X,{\F}}(x)$ для всех $x\in\cN^0$ и поэтому ${\nu}_{X,{\F}}$ --- корректно определенная полунорма на $\cN^0$.
\end{rk}

Ясно, что $\mathbb {OSN}(\cN,\cN^0)=\mathbb{SN}(M(\cN),\cN^0)$ и $\mathbb{OA}(\cN,\cN^0)=\mathbb{A}(M(\cN),\cN^0)$.

\begin{dfn}\label{dfn:pseme}
Для $(X,\F)\in \mathbb{OA}(\cN,\cN^0)$,
рассмотрим функцию
$d_{X,\F}:\cN^0\times \cN^0\to [0,+\infty)$
$$
d_{X,\F}(x,y)^2:=\nu_{X,\F}(x-y)^2=
\sup_k 
\sum_{i=k}^\infty |\f_k\left(\<x-y,x_i\>\right)|^2,\quad x,y\in \cN^0. 
$$
Будем писать $d_{X,\F}\in \mathbb{OPM}(\cN,\cN^0)$ (=$\mathbb{PM}(M(\cN),\cN^0)$).
\end{dfn}

Очевидно, $d_{X,\F}$ являются \emph{полуметриками} в смысле  \cite[определение 2.10]{Troitsky2020JMAA} (и \cite[
Глава IX, \S 1]{BourbakiTop2}), следовательно, они образуют равномерную структуру.

Определение \emph{вполне ограниченных} множеств 
относительно рассматриваемой равномерной структуры
(или для системы $\mathbb{OPM}(\cN,\cN^0)$) принимает следующий вид.

\begin{dfn}\label{dfn:totbaundset}
\rm
Множество $Y\subseteq \cN^0 \subseteq  \cN(\subseteq  M(\cN))$ называется \emph{вполне ограниченным}
по отношению к этой равномерной структуре, если для любого $(X,\F)$, 
где $X \subseteq M(\cN)$ является внешней $\cN^0$-допустимой,
и любого $\e>0$ найдется конечный набор
$y_1,\dots,y_n$
элементов $Y$, такой что множества
$$
\left\{ y\in Y\,|\, d_{X,\F}(y_i,y)<\e\right\}
$$  
образуют покрытие $Y$. Этот конечный набор называется
$\e$\emph{-сетью в $Y$ для} $d_{X,\F}$.

Мы будем говорить в этом случае, что $Y$ является \emph{внешне} $(\cN,\cN^0)$-\emph{вполне ограниченным}.
Ясно, что это эквивалентно $(M(\cN),\cN^0)$-вполне ограниченности $Y$.
\end{dfn}

Очевидно, $(M(\cN),\cN^0)$-вполне ограниченное множество является $(\cN,\cN^0)$-вполне \linebreak
 ограниченным, поэтому все результаты, устанавливающие, что $(\cN,\cN)$-вполне ограниченность образа единичного шара влечет $\A$-компактность соответствующего оператора, допускающего сопряженный (\cite[Теорема 2.13]{Troitsky2020JMAA}, \cite[Теорема 3.5]{TroitFuf2020})
остаются верными, если рассматривать
$(M(\cN),\cN)$-вполне ограниченность. Наша нынешняя цель состоит в том, чтобы показать, что существует пример, для которого данный вывод неверен --- а именно, существует (коммутативная) $C^*$-алгебра $\A$, такая, что если её рассматривать как модуль над собой, то тождественный оператор $Id:\A\to\A$ не является $\A$-компактным, но единичный шар в $\A$ (и, поэтому, образ единичного шара под действием тождественного оператора) является $(\cN,\cN)$-вполне ограниченным. Более того, мы покажем, что он даже $(M(\cN),\cN)$-вполне ограничен.

Наконец, напомним понятие стандартного фрейма в гильбертовом $C^*$-модуле (см., например,  \cite{FrankLarson1999}, \cite{FrankLarson2002}).

\begin{dfn}\label{dfn:fr}
Пусть $\cN$ --- гильбертов $C^*$-модуль над унитальной $C^*$-алгеброй $\A$ и $J$ --- некоторое множество. Семейство $\{x_j\}_{j\in J}$ элементов $\cN$ называется стандартным фреймом в $\cN$, если найдутся такие положительные постоянные $c_1, c_2$, что
для любого $x\in\cN$
ряд
$\sum\limits_j \<x,x_j\>\<x_j,x\>$
сходится по норме в $\A$
и справедливо следующее неравенство:
 $$c_1\<x,x\>\le \sum\limits_j \<x,x_j\>\<x_j,x\>\le c_2\<x,x\>.$$
\end{dfn}

\section{$(M(\cN),\cN)$-вполне ограниченность не влечет $\A$-компактность}\label{sec:TotBoundNoAcomp} 

\begin{lem}\label{cor:noncgen}
Как модуль над собой, $C_0(K)$ не является счетнопорожденным.
\end{lem}

\begin{proof}
Действительно, если $\{f_n\}_{n\in\mathbb N}$ --- счетная система образующих для $C_0(K)$ и компакт $K_n$ --- носитель $f_n$, то всякая функция из $C_0(K)$
обращается в ноль вне outside $K'$, где $K'$ --- компакт, содержащий все $K_n$,
но это неверно, потому что найдутся $z\in K\setminus K'$ и, так как $K$ является пространством Урысона, функция $f\in C_0(K)$ такая, что $f(z)=1$. 
Более точно, такая функция найдется, потому что функции из $C_0(K)$ --- это в точности непрерывные функции на $\widetilde{K}=K\cup\{t_\infty\}$ --- одноточечной компактификации $K$,
такие что $f(t_\infty)=0$, поэтому свойство Урысона дает функцию $f\in C(\widetilde{K})$ такую, что $f(z)=1$, $f(t_\infty)=0$.
\end{proof}

\begin{cor}\label{cor:idncomp}
Тождественный оператор $Id:C_0(K)\to C_0(K)$ не является $C_0(K)$-компактным.
\end{cor}

\begin{proof}
Действительно, как было показано в 
 \cite[лемма 1.10]{Troitsky2020JMAA},
для $C^*$-алгебры $\A$ образ $\A$-компактного оператора между двумя гильбертовыми $C^*$-модулями над $\A$ должен быть счетнопорожденным.
\end{proof}

\begin{lem}\label{lem:compsupm}
Если $\mu$ --- счетно-аддитивная положительная конечная Радоновская мера на $K$ (т.е. положительный непрерывный линейный функционал на $C_0(K)$, см. \cite[436K]{Frem4})
и $K'$ --- ее носитель, т.е. 
$\mu(A)=0$ для любого борелевского множества $A\subset K\setminus K'$ и, следовательно,
$\mu(f)=\int\limits_Kf(t)d\mu(t)=0$ для всех $f\in C_0(K)$ таких, что $f(t)=0$ при $t\in K'$, то $K'$ --- компакт.
\end{lem}

\begin{proof}
Напомним, что в силу \cite[411B, 411H]{Frem4} для любого борелевского множества $E$ выполнено $\mu(E)=\sup\{\mu(D), D\subset E, D\,\, \text{компактно}\}$.

Предположим, что $K'$ --- не компакт. Тогда он не может быть покрыт счетным множеством компактов, так как если $K'\subset\bigcup_{n=1}^\infty K_n$, где $K_n$ --- компакты, то носитель $\mu$ содержится в компактном множестве, содержащем все $K_n$.

Поэтому найдется несчетное семейство компактных множеств $K_\a$, упорядоченных по строгому включению, таких что $\mu(K_\beta\setminus K_\a)>0$ при $K_\a\subset K_\beta$. На самом деле, это семейство может быть занумеровано порядковыми числами, меньшими, чем $\om_1$: выберем произвольный компакт $K_1$ с $\mu(K_1)>0$, и если мы уже нашли $K_\a$ для всех $\a<\beta$, то $\overline{\bigcup_{\a<\beta} K_\a}$ --- компактное множество, так как это замыкание не более чем счетного объединения компактов, и, так как $\mu$ не имеет компактного носителя, найдется компакт $K_\beta$, $\overline{\bigcup_{\a<\beta} K_\a}\subset K_\beta$, такой что $\mu( K_\beta\setminus \overline{\bigcup_{\a<\beta} K_\a})>0$.
Действительно, $K\setminus \overline{\bigcup_{\a<\beta} K_\a}$ --- борелевское множество и $\mu( K\setminus \overline{\bigcup_{\a<\beta} K_\a})>0$, поэтому
найдется компакт $D_\beta\subset K\setminus \overline{\bigcup_{\a<\beta} K_\a}$ с $\mu(D_\beta)>0$.
Можем положить $K_\beta=\overline{\bigcup_{\a<\beta} K_\a}\cup D_\beta$.

Для каждого порядкового числа $\a<\om_1$ существует следующее за ним порядковое число $\a+1$, поэтому мы можем ввести несчетное множество непересекающихся борелевских множеств $Y_\a=K_{\a+1}\setminus K_{\a}$, таких что $0<\mu(Y_\a)<+\infty$.

Если для всякого $n\in\mathbb N$ найдется только конечный набор множеств $Y_\a$ таких, что $\mu(Y_\a)\ge\frac{1}{n}$, то существует только не более чем счетное число $Y_\a$ с $\mu(Y_\a)\ne0$, но это неверно, поэтому найдется $k\in\mathbb N$ такое, что для бесконечного набора (даже несчетного) борелевских множеств $Y_\a$ выполнено $\mu(Y_\a)\ge \frac{1}{k}$. Возьмем последовательность $Y_{\a_j}$, $j\in\mathbb N$, $\a_j<\a_{j+1}$, таких множеств. В силу того, что $\mu$ --- Радоновская мера и $Y_{\a_j}$ --- борелевские множества,
найдутся компактные множества $F_{\a_j}\subset Y_{\a_j}$ такие, что $\mu(F_{\a_j})\ge \frac{1}{2k}$.

$F_{\a_j}\subset K_{\a_j+1}$, $K_{\a_j+1}$ --- компакт. Из \cite[3.3.2 и 3.3.3]{Engel} следует, что найдется неотрицательная функция $g_j\in C_0(K)$ такая, что $g_j(t)=1$ при $t\in K_{\a_j+1}$. Очевидно, $||g_j||_{C_0(K)}=1$. Но $||\mu||_{C_0(K)^*}\ge|\mu(g_j)|=\int\limits_Kg_j(t)d\mu(t)\ge\int_{\bigsqcup\limits_{l=1}^jF_{\a_l}}g_j(t)d\mu(t)\ge \frac{j}{2k}\to\infty$ при $j\to\infty$, поэтому $\mu$ не является непрерывным линейным функционалом на $C_0(K)$. Противоречие, следовательно, $\mu$ имеет компактный носитель.
\end{proof}

\begin{rk}
Можно также получить противоречие в предыдущем доказательстве, если заметить, что все компакты $K_{\a_j}$ содержаться в некотором компакте $K'$, поэтому $\mu(K')\ge\mu(\bigsqcup\limits_{l=1}^\infty F_{\a_l})=\sum\limits_{l=1}^\infty \mu(F_{\a_l})=\sum\limits_{l=1}^\infty  \frac{1}{2k}=+\infty$, но $\mu(K')$ должно быть конечной величиной.
\end{rk}

Возьмем полунорму
$
\nu_{X,\F}(x)^2:=\sup_k 
\sum_{i=k}^\infty |\f_k\left(\<x,x_i\>\right)|^2
$
с  $X=\{x_i\}$ --- допустимой системой элементов мультипликаторного модуля (то есть,
внешней
допустимой системой). Так как в качестве модуля мы берем алгебру $C_0(K)$ с внутренним произведением $\<a,b\>=a^*b$, мультипликаторный модуль совпадает с пространством $C_b(K)$ всех ограниченных непрерывных функций на $K$.

\begin{teo}
Единичный шар в $C_0(K)$ (и, следовательно, образ единичного шара под действием тождественного оператора $Id:C_0(K)\to C_0(K)$) вполне ограничен относительно полунормы $\nu_{X,\F}$ для любой $(C_b(K),C_0(K))$-допустимой системы $X$ и любой последовательности состояний $\F$. Иными словами, он $(C_b(K),C_0(K))$-вполне ограничен. 
\end{teo}

\begin{proof}
Для любого $k\in\mathbb N$, так как $\f_k$ является непрерывным положительным линейным функционалом на $C_0(K)$, найдется счетно-аддитивная положительная конечная Радоновская мера $\mu_k$ на $K$ с компактным носителем $K_k$, и, следовательно, найдется компакт
$K'$,
который содержит носители всех мер $\mu_k$. 

Поэтому
$
\nu_{X,\F}(x)^2=
\sup_k 
\sum_{i=k}^\infty |
\int\limits_K\overline{x(t)}\cdot x_i(t)d\mu_k(t)
|^2
=
\sup_k 
\sum_{i=k}^\infty |
\int\limits_{K'}\overline{x(t)}\cdot x_i(t)d\mu_k(t)
|^2
$.
Эта величина, в действительности, является полунормой на модуле $C(K')$
в ``старом'' смысле (т.е. в смысле \cite[определения 2.1 и 2.4]{Troitsky2020JMAA}, без рассмотрения мультипликаторного модуля), так как сужения $x$, $x_i$ на $K'$ являются непрерывными функциями, $\f_i$ являются непрерывными линейными функционалами на $C(K')$ и система $\{x_i\}$ (более точно, система сужений функций) является $C(K')$-допустимой системой элементов из $C(K')$.
Действительно, условия

\begin{enumerate}
\item[1')] 
при каждом $x\in
 C(K')
 $
ряд $\sum_i \<x,x_i\>\<x_i,x\>$ сходится по норме;
\item[2')]
его сумма ограничена $\<x,x\>$; 
\item[3')]$\|x_i\|\le 1$ для любого $i$.
\end{enumerate}
в $C(K')$ выполнены, так как каждая функция $x\in C(K')$ имеет продолжение $\dot{x}\in C_0(K)$ и условия 

\begin{enumerate}
\item[1)] 
при каждом $\dot{x}\in
 C_0(K)
 $,
ряд $\sum_i \<\dot{x},x_i\>\<x_i,\dot{x}\>$ сходится по норме;
\item[2)]
его сумма ограничена $\<\dot{x},\dot{x}\>$; 
\item[3)]$\|x_i\|\le 1$ для любого $i$.
\end{enumerate}
в $C_0(K)$ выполнены, так как $\{x_i\}$ является допустимой системой. Поэтому, 1') следует из 1) так как равномерно сходящийся ряд также сходится равномерно на подмножествах, 2') следует из 2) так как это условие --- поточечное, и если оно выполнено на множества, то оно также выполнено и на подмножествах, 3') следует из 3) очевидно.

Тождественный оператор $Id:C(K')\to C(K')$ является $C(K')$-компактным, потому что эта $C^*$-алгебра --- унитальна. Также, она счетнопорождена как $C(K')$-модуль. Поэтому, по \cite[теорема 2.13]{Troitsky2020JMAA}, единичный шар $C(K')$ вполне ограничен относительно $\nu_{X,\F}$, то есть, для любого $\e>0$ найдется конечная $\e$-сеть
$\{y_1,\dots,y_m\}$ из элементов $C(K')$ для $\nu_{X,\F}$. 
По следствию \ref{extrest},
найдется набор элементов
$\{\dot{y_1},\dots,\dot{y_m}\}$ из $C_0(K)$, которые являются продолжениями $\{y_1,\dots,y_m\}$. В силу того, что $\nu_{X,\F}$ на $C_0(K)$ в действительности вычисляется как полунорма на $C(K')$, система $\{\dot{y_1},\dots,\dot{y_m}\}$ является $\e$-сетью для единичного шара $C_0(K)$.
\end{proof}

\begin{rk}
Отметим, что если для $\nu_{X,\F}$ мы рассмотрим только допустимые системы, то есть, в нашем случае, функции из $C_0(K)$, то для всякого счетного набора таких функций найдется компакт $K'$ такой, что они обращаются в ноль вне $K'$, поэтому $\nu_{X,\F}$ в действительности вычисляется как полунорма $C(K')$ без применения свойств Радоновских мер на $K$, т.е. без леммы \ref{lem:compsupm}.
В случае внешних допустимых систем, т.е. в случае $C_b(K)$, этот аргумент использовать нельзя.
\end{rk}

\section{Отсутствие фреймов}\label{sec:nofr}

Очевидно, $C^*$-алгебру $\A$ (если она не унитальна) можно рассматривать как гильбертов $C^*$-модуль над ее унитализацией $\dot{\A}$.

\begin{teo}\label{nofr}
В $\dot{C_0}(K)$-модуле $C_0(K)$ нет стандартного фрейма.
\end{teo}

\begin{proof}
Предположим, что существует фрейм $\{x_j\}_{j\in J}$ из элементов $C_0(K)$.

Возьмем произвольную точку $z_1\in K$. Найдется функция $g_1\in C_0(K)$ такая, что $g_1(z_1)=1$. Существует непустое не более чем счетное множество $\{x_j\}_{j\in J_1}$ элементов фрейма, такое что $x_j(z_1)\ne0$ в силу того, что

$$c_1\<g_1,g_1\>\le \sum\limits_j \<g_1,x_j\>\<x_j,g_1\>\le c_2\<g_1,g_1\>$$
 
и, следовательно,

$$c_1|g_1(z_1)|^2\le \sum\limits_j |g_1(z_1)|^2 |x_j(z_1)|^2\le c_2|g_1(z_1)|^2$$

и
 
$$c_1\le \sum\limits_j |x_j(z_1)|^2\le c_2.
$$
  
Так как
$\{x_j\}_{j\in J_1}$ --- не более чем счетное множество функций из $C_0(K)$, найдется компакт $K_1$ такой, что все эти функции обращаются в ноль вне $K_1$.

Применим метод математической индукции: пусть мы уже нашли точки $z_l\in K$, не более чем счетные множества индексов $J_l\subset J$, компактные множества $K_l$, $l=1,\dots,m$, такие что $K_1\subset K_2\subset\dots\subset K_m$, $x_j$ обращаются в ноль вне $K_m$ для всех $j\in J_1\cup\dots\cup J_m$.

Возьмем произвольную $z_{m+1}\in K\setminus K_m$.
Найдется функция $g_{m+1}\in C_0(K)$ такая, что $g_{m+1}(z_{m+1})=1$. Существует непустое не более чем счетное множество $\{x_j\}_{j\in J_{m+1}}$ элементов фрейма, такое что $x_j(z_{m+1})\ne0$ в силу того, что

$$c_1\<g_{m+1},g_{m+1}\>\le \sum\limits_j \<g_{m+1},x_j\>\<x_j,g_{m+1}\>\le c_2\<g_{m+1},g_{m+1}\>$$

и, следовательно,
 
$$c_1|g_{m+1}(z_{m+1})|^2\le \sum\limits_j |g_{m+1}(z_{m+1})|^2 |x_j(z_{m+1})|^2\le c_2|g_{m+1}(z_{m+1})|^2$$

и 
 
$$c_1\le \sum\limits_j |x_j(z_{m+1})|^2\le c_2.
$$
 
Отметим, что $J_{m+1}\cap J_l=\o$ для всех $l=1,\dots,m$.
 
Так как
$\{x_j\}_{j\in J_1\cup\dots\cup J_{m+1}}$ --- не более чем счетное множество функций из $C_0(K)$, найдется компакт $K_{m+1}$ такой, что все эти функции обращаются в ноль вне $K_{m+1}$ и $K_m\subset K_{m+1}$.
И так далее.

Таким образом, построены последовательности $\{z_l\}$, $\{J_l\}$, $\{K_l\}$. Найдется компакт $K'$ такой, что $\{z_l\}\subset K'$ и найдется такая функция $g\in C_0(K)$ что $g(t)=1$ при $t\in K'$. Так как $\{x_j\}$ --- фрейм, ряд $\sum\limits_j \<g,x_j\>\<x_j,g\>$ сходится в $ C_0(K)$ (т.е. равномерно). Но всякая частичная сумма $S_N$ этого ряда состоит из слагаемых только из конечного набора индексных множеств $J_1\cup J_2\cup\dots\cup J_{m-1}$, поэтому $S_N$ обращается в ноль вне $K_{m-1}$ и не равна нулю только для конечного множества точек из $\{z_l\}$: $z_1,\dots,z_{m-1}$.
Для всей суммы в точке $z_m$ имеем
$\sum\limits_j \<g,x_j\>\<x_j,g\>(z_m)=\sum\limits_j |g(z_{m})|^2 |x_j(z_{m})|^2\ge c_1|g(z_{m})|^2=c_1>0$,
поэтому ряд не сходится равномерно. Противоречие.
\end{proof}

\begin{cor}
Из теоремы \ref{nofr}
и результатов Франка и Ларсона (\cite[3.5, 4.1 и 5.3]{FrankLarson2002}, см. также \cite[Теорема 1.1]{HLi2010}) следует, что $C_0(K)$ не может быть представлен как прямое слагаемое стандартного модуля $\bigoplus\limits_{\a\in\Lambda}\dot{C_0}(K)$ никакой мощности.
\end{cor}

\begin{rk}
Несмотря на то, что $C_0(K)$ не имеет стандартных фреймов, он, тем не менее, имеет внешние фреймы (т.е. фреймы, чьи элементы могут лежать в $M(\cN)$, не только в $\cN$, подробнее см. \cite{ArBak2017}). Действительно, для любой $C^*$-алгебры $\A$, рассматриваемой как модуль над собой, можно взять, например, внешний фрейм, состоящий лишь из единичного элемента алгебры $M(\A)$.
\end{rk}

\end{document}